\documentclass[11pt]{article}
\usepackage{amsfonts}
\usepackage{bbm}
\usepackage{mathrsfs}
\usepackage{amscd}
\usepackage{color}
\usepackage{amsmath,amsfonts,amssymb,amscd}
\usepackage{indentfirst,graphics,epsfig,psfrag}
\input{epsf}
\usepackage{ifpdf}
\usepackage{enumerate}
\usepackage{appendix}
\usepackage{enumerate}

\usepackage{lineno}

\makeatletter \@addtoreset{figure}{section} \makeatother
\makeatletter
\long\def\@makecaption#1#2{%
   \vskip 10\p@
   \setbox\@tempboxa\hbox{{#1}\ \ #2}%
   \ifdim \wd\@tempboxa >\hsize

       {#1}\ \ #2\par
   \else
       \hbox to\hsize{\hfil\box\@tempboxa\hfil}%
   \fi}
\makeatother

\begin{document}
\title{\textbf{Results on resistance distance and Kirchhoff index
of graphs with generalized pockets }}
\author{
\small  Qun Liu~$^{a,b}$\\
\setcounter{footnote}{-1 }\thanks{Corresponding
author: Zhongzhi Zhang, E-mail address: liuqun09@yeah.net, zhangzz@fudan.edu.cn}\\
\small  a. School of Computer Science, Fudan University, Shanghai 200433, China\\
\small  b. School of Mathematics and Statistics, 
Hexi University, Gansu, Zhangye, 734000, P.R. China\\}
\date{}
\maketitle
\begin{abstract}
Let $F,H_{v}$ be simple connected graphs on $n$ and $m + 1$ vertices, respectively. Let $v$ be a specified vertex of $H_{v}$ and $u_{1},u_{2},\ldots u_{v}\in F$. Then the graph $G= G[F,u_{1},...,u_{k}, H_{v}]$ obtained by taking one copy of $F$ and $k$ copies of $H_{v}$, and then attaching the $i$th copy of $H_{v}$ to the vertex $u_{i}, i = 1,...,k$, at the vertex $v$ of $H_{v}$(identify $u_{i}$ with the vertex $v$ of the $i$th copy) is called a graph with $k$ pockets. In $\cite{SB2018}$, Barik gave the Laplacian spectrum for more general cases. In this paper, we derive closed-form formulas for resistance distance
and Kirchhoff index of $G= G[F,u_{1},...,u_{k}, H_{v}]$ in terms of the resistance distance and Kirchhoff index
$F$ and $H_{uv}$, respectively.
\\[2mm]
\\[1mm]
\end{abstract}

\section{ Introduction }
All graphs considered in this paper are simple and undirected.
The resistance distance between vertices $u$ and $v$ of $G$ was defined
by Klein and Randi$\acute{c}$ $\cite{KR}$ to be the effective resistance between nodes $u$ and $v$ as computed with Ohm's law when all the edges of $G$
are considered to be unit resistors.
The Kirchhoff index $Kf(G)$ was defined in $\cite{KR}$
as $Kf(G)=\sum_{u<v}r_{uv}$, where $r_{uv}(G)$ denote the resistance distance between $u$ and $v$ in $G$. Resistance distance are, in
fact, intrinsic to the graph, with some nice purely mathematical
interpretations and other interpretations. The Kirchhoff
index was introduced in chemistry as a better alternative
to other parameters used for discriminating different molecules
with similar shapes and structures $\cite{KR}$. The resistance distance and the Kirchhoff index have attracted extensive attention due to its wide applications in physics, chemistry and others.
Up till now, many results on the resistance distance
and the Kirchhoff index are obtained. See
$(\cite{JZ}, \cite{JB1}, \cite{JB2}, \cite{SWZhB})$ and the references therein
to know more. However, the resistance
distance and Kirchhoff index of the graph is, in general, a difficult
thing from the computational point of view. Therefore, the bigger is
the graph, the more difficult is to compute the resistance distance
and Kirchhoff index, so a common strategy is to consider complex graph
as composite graph, and to find relations between the resistance distance
and Kirchhoff index of the original graphs.


Let $G=(V(G),E(G))$ be a graph with vertex set $V(G)$ and edge set $E(G)$. Let $d_{i}$ be the degree of vertex $i$ in $G$ and $D_{G}=diag(d_{1}, d_{2},\cdots d_{|V(G)|})$ the diagonal matrix with all vertex degrees of $G$ as its diagonal entries. For a graph $G$, let $A_{G}$ and $B_{G}$ denote the adjacency matrix and vertex-edge incidence matrix of $G$, respectively. The matrix $L_{G}=D_{G}-A_{G}$ is called the Laplacian matrix of $G$, where $D_{G}$ is the diagonal matrix of vertex degrees of $G$. We use $\mu_{1}(G)\geq u_{2}(G)\geq\cdots\geq\mu_{n}(G)=0$ to denote the eigenvalues of $L_{G}$. For other undefined notations and terminology from graph theory, the readers may refer to $\cite{R. B. Bapat}$ and the references therein $(\cite{ChZh}-\cite{YK})$.

The computation of resistance distance between two nodes in a resistor network is a classical problem in electric theory and graph theory. For certain families of graphs it is possible to identify a graph by looking at the resistance distance and Kirchhoff index. More generally, this is not possible. In some cases, the
resistance distance and Kirchhoff index of a relatively larger graph can be described in terms of the resistance distance and Kirchhoff index of some smaller(and simpler) graphs using some simple graph operations. There are results that discuss the resistance distance and Kirchhoff index of graphs obtained by means of some operations on graphs like join, graph products, corona and many variants of corona(like edge corona, neighborhood corona, edge, neighborhood corona,etc.). For such operations often it is possible to describe the resistance distance and Kirchhoff index of the resulting graph using the resistance distance and Kirchhoff index of the corresponding constituting graph, see  $(\cite{BYZhZH},\cite{LZhB})$ for
reference. This paper consider the resistance distance and Kirchhoff index of the graph operations below, which come
from \cite{SB}.

{\bf{Definition 1 \cite{SB}}}  Let $F$, $H_{v}$ be connected graphs, $v$ be a specified vertex of $H_{v}$ and $u_{1},u_{2},...,u_{k}\in F$.
Let $G=G[F,u_{1},u_{2},...,u_{k},H_{v}]$ be the graph obtained by taking one copy of $F$ and $k$ copies of $H_{v}$, and then
attaching the $i$th copy of $H_{v}$ to the vertex $u_{i}$, $i=1,2,...,k$, at the vertex $v$ of $H_{v}$(identify $u_{i}$ with the vertex$v$ of the $i$th copy). Then the copies of the graph $H_{v}$ that are attached to the vertices $u_{i}$, $i=1,2,...,k$ are referred to as pockets,
and $G$ is described as a graph with $k$ pockets.

Barik $\cite{SB}$ has described the $L$-spectrum of
$G=G[F,u_{1},u_{2},...,u_{k},H_{v}]$ using the
$L$-spectrum of $F$ and $H_{v}$ in a particular
case when $deg(v)=m$. Recently Barik and Sahoo
$\cite{SB2018}$ have described the Laplacian spectrum of more
such graphs relaxing condition $deg(v)=m$. Let
$deg(v)=l$, $1\leq l\leq m$. In this case, we denote
$G=G[F,u_{1},u_{2},...,u_{k},H_{v}]$ more precisely
by $G=G[F,u_{1},u_{2},...,u_{k};H_{v},l]$. When
$k=n$, we denote simply by $G=G[F;H_{v},l]$. If
$deg(v)=l$, $1\leq l\leq m$, let $N(v)=\{
v_{1},v_{2},\ldots,v_{l}\}\subset V(H_{v})$ be the
neighbourhood set of $v$ in $H_{v}$. Let $H_{1}$ be
the subgraph of $H(v)$ induced by the vertices in
$N(v)$ and $H_{2}$ be the subgraph of $H_{v}$ induced
by the vertices which are in $V(H_{v})\backslash(N(v)\cup\{v\})$. When $H_{v}=H_{1}\vee(H_{2}+\{v\})$, we describe the
resistance distance and Kirchhoff index of
$G=G[F,u_{1},u_{2},...,u_{k},H_{v}]$. The results are contained in Section 3 of the article.

Further, when $F=F_{1}\vee F_{2}$,  where
$F_{1}$ is the subgraph of $F$ induced by the vertices
$u_{1},u_{2},\ldots,u_{k}$ and $F_{2}$ is the subgraph of $F$
induced by the vertices $u_{k+1},u_{k+2},\ldots,u_{n}$. In this case, we describe the resistance distance and Kirchhoff index of $G[F,u_{1},u_{2},...,u_{k};H_{v},l]$. These results are contained in Section 4.

\section{Preliminaries}
The $\{1\}$-inverse of $M$ is a matrix $X$ such that $MXM=M$. If $M$ is singular, then it has infinite
$\{1\}$-inverse $\cite{A. Ben-Israel}$.
For a square matrix $M$, the group inverse of $M$, denoted by $M^{\#}$, is the unique matrix $X$ such that $MXM=M$, $XMX=X$ and $MX=XM$. It is known that $M^{\#}$ exists if and only if $rank(M)=rank(M^{2})$ $(\cite{A. Ben-Israel},\cite{BSZHW})$. If $M$ is real symmetric, then $M^{\#}$ exists and $M^{\#}$ is a symmetric
$\{1\}$-inverse of $M$. Actually, $M^{\#}$ is equal to the Moore-Penrose inverse of $M$ since $M$ is symmetric $\cite{BSZHW}$.

It is known that resistance distances in a connected graph $G$ can be obtained from any $\{1\}$- inverse of $G$ $(\cite{BG})$. We use $M^{(1)}$ to denote any $\{{1}\}$-inverse of a matrix $M$, and let $(M)_{uv}$ denote the $(u,v)$-entry of $M$.

{\bf Lemma 2.1} $(\cite{BSZHW})$\ Let $G$ be a connected graph. Then
$$
r_{uv}(G)=(L^{(1)}_{G})_{uu}+(L^{(1)}_{G})_{vv}-(L^{(1)}_{G})_{uv}
-(L^{(1)}_{G})_{vu}=(L^{\#}_{G})_{uu}+(L^{\#}_{G})_{vv}-2(L^{\#}_{G})_{uv}.$$

Let $1_{n}$ denote the column vector of dimension $n$ with all the entries equal one. We will often use $1$ to denote all-ones column vector if the dimension can be read from the context.

{\bf Lemma 2.2} $(\cite{BYZhZH})$ \ For any graph, we have
$L^{\#}_{G}1=0.$

{\bf Lemma 2.3} $(\cite{FZZh})$ \ Let
 \[
\begin{array}{crl}
M=\Bigg(
  \begin{array}{cccccccccccccccc}
   A& B  \\
   C & D \\
 \end{array}
  \Bigg)
\end{array}
\]
be a nonsingular matrix. If $A$ and $D$ are nonsingular,
then
\[
\begin{array}{crl}
M^{-1}&=&\Bigg(
  \begin{array}{cccccccccccccccc}
   A^{-1}+A^{-1}BS^{-1}CA^{-1}&-A^{-1}BS^{-1} \\
   -S^{-1}CA^{-1} & S^{-1}\\
 \end{array}
  \Bigg)
\\&=&\Bigg(
  \begin{array}{cccccccccccccccc}
   (A-BD^{-1}C)^{-1} &-A^{-1}BS^{-1} \\
   -S^{-1}CA^{-1} & S^{-1}\\
 \end{array}
  \Bigg),
\end{array}
\]
where $S=D-CA^{-1}B.$




{\bf Lemma 2.4} $(\cite{LZhB})$ \ Let $L$ be the Laplacian
matrix a graph of order $n$. For any $a>0$, we have
$$(L+aI_{n}-\frac{a}{n}j_{n\times n})^{\#}
=(L+aI)^{-1}-\frac{1}{an}j_{n\times n}.$$

{\bf Lemma 2.5} $(\cite{SWZhB})$ \ Let $G$ be a connected graph on $n$ vertices. Then
$$Kf(G)=ntr(L^{(1)}_{G})-1^{T}L^{(1)}_{G}1=ntr(L^{\#}_{G}).$$


{\bf Lemma 2.6} $(\cite{QL})$ \ Let
 \[
\begin{array}{crl}
L=\left(
  \begin{array}{cccccccccccccccc}
   A& B \\
   B^{T} & D \\
 \end{array}
  \right)
\end{array}
\]
be the Laplacian matrix of a connected graph.
If $D$ is nonsingular, then
 \[
\begin{array}{crl}
X=\left(
  \begin{array}{cccccccccccccccc}
   H^{\#}& -H^{\#}BD^{-1} \\
   -D^{-1}B^{T}H^{\#} & D^{-1}+D^{-1}B^{T}H^{\#}BD^{-1} \\
 \end{array}
  \right)
\end{array}
\]
is a symmetric $\{1\}$-inverse of $L$,
where $H=A-BD^{-1}B^{T}$.

\section{The resistance distance and Kirchhoff index of $G[F;H_{v},l]$}
Let $F$ be a connected graph with vertex set $\{u_{1},u_{2},...,u_{n}\}$. Let $H_{v}$ be a
connected graph on $m+1$ vertices with a specified
vertex $v$ and $V(H_{v})=\{v_{1},v_{2},...,v_{m},v\}$.
Let $G=G[F;H_{v},l]$. Note that $G$ has $n(m+1)$ vertices. Let $deg(v)=l$, $1\leq l\leq m$. With loss of generality
assume that $N(v)=\{v_{1},v_{2},...,v_{l}\}$.
Let $H_{1}$ be
the subgraph of $H_{v}$ induced by the vertices in
$\{v_{1},v_{2},...,v_{l}\}$ and $H_{2}$ be the subgraph of $H_{v}$ induced
by the vertices $\{v_{l+1},v_{l+2},...,v_{m}\}$.
Suppose that $H_{v}=H_{1}\vee(H_{2}+\{v\})$.
In this section, we focus on determing the resistance distance and Kirchhoff index of $G[F;H_{v},l]$ in terms of the resistance distance and Kirchhoff index of $F, H_{1}$ and $H_{2}$.

{\bf Theorem 3.1} \  Let $G[F;H_{v},l]$
be the graph as described above. Suppose that $H_{v}=H_{1}\vee(H_{2}+\{v\})$. Let the Laplacian spectrum of
$H_{1}$ and $H_{2}$ be
$\sigma(H_{1})=(0=\mu_{1},\mu_{2},\ldots \mu_{l})$ and
$\sigma(H_{2})=(0=\nu_{1},\nu_{2},\ldots \nu_{m-l})$. Then
$G[F;H_{v},l]$ have the resistance distance and Kirchhoff index
as follows:

(i) For any $i,j\in V(F)$, we have
\begin{eqnarray*}
r_{ij}(G[F;H_{v},l])=(L^{\#}(F))_{ii}+(L^{\#}(F))_{jj}
-2(L^{\#}(F))_{ij}.
\end{eqnarray*}

(ii) For any $i\in V(F)$, $j\in V(H_{1})$, we have
\begin{eqnarray*}
r_{ij}(G[F;H_{v},l])&=&(L^{\#}(F))_{ii}+
[(L(H_{1})+(m-l+1)I_{l}-\frac{m-l}{l}J_{l\times l})^{-1}\otimes I_{n}]_{jj}-2(L^{\#}(F))_{ij}.
\end{eqnarray*}

(iii) For any $i\in V(F)$, $j\in V(H_{2})$, we have
\begin{eqnarray*}
r_{ij}(G[F;H_{v},l])&=&(L^{\#}(F))_{ii}+
[(L(H_{2})+lI_{m-l}-\frac{l}{m-l+1}J_{(m-l)\times (m-l)})^{-1}\otimes I_{n}]_{jj}-2(L^{\#}(F))_{ij}.
\end{eqnarray*}

(iv)For any $i\in V(H_{1})$, $j\in V(H_{2})$, we have

$r_{ij}(G[F;H_{v},l])$
\begin{eqnarray*}
&=&[(L(H_{1})+(m-l+1)I_{l}-\frac{m-l}{l}J_{l\times l})^{-1}\otimes I_{n}]_{ii}+\left[(L(H_{2})+lI_{m-l}\right.\\&&
\left.-\frac{l}{m-l+1}J_{(m-l)\times (m-l)})^{-1}\otimes I_{n}]_{jj}-2[(L(H_{1})+(m-l+1)I_{l}-\frac{m-l}{l}J_{l\times l})^{-1}\otimes I_{n}\right]_{ij}.
\end{eqnarray*}

(v)For any $i\in V(H_{2})$, $j\in V(H_{1})$, we have
\begin{eqnarray*}
r_{ij}(G[F;H_{v},l])&=&[(L(H_{2})+lI_{m-l}-\frac{l}{m-l+1}J_{(m-l)\times (m-l)})^{-1}\otimes I_{n}]_{ii}+\left[(L(H_{1})+(m-l+1)I_{l}\right.\\&&
\left.-\frac{m-l}{l}J_{l\times l})^{-1}\otimes I_{n}]_{jj}-2[(L(H_{2})+lI_{m-l}-\frac{l}{m-l+1}J_{(m-l)\times (m-l)})^{-1}\otimes I_{n}\right]_{ij}.
\end{eqnarray*}

(iv)
$Kf(G[F;H_{v},l])$
\begin{eqnarray*}
&=&n(m+1)\left(\frac{m+1}{n}Kf(F)
+(n\sum_{i=2}^{l}\frac{1}{\mu_{i}(H_{1})+(m-l+1)}+n)+
(n\sum_{i=2}^{m-l}\frac{1}{\nu_{i}(H_{2})+l}\right.\\&&
\left.+\frac{nl}{m-l+1})\right)-\left((m-l)^{2}\frac{m-l+1}{l}+l^{2}
\right).
\end{eqnarray*}

{\bf Proof} \ Let $v^{i}_{j}$ denote the $j$th vertex of $H_{v}$
in the $i$th copy of $H_{v}$ in $G$, for $i=1,2,\ldots,n$;
$j=1,2,\ldots,m$, and let $V_{j}(H_{v})=\{v^{1}_{j},v^{2}_{j},
\ldots,v^{n}_{j}\}$. Then $V(F)\cup (\cup_{j=1}^{m}V_{j}(H_{v}))$
is a partition of $V(G)$. Using this partition, the Laplacian
matrix of $G=G[F;H_{v},l]$ can be expressed as
\[
\begin{array}{crl}L(G)=
\left(
  \begin{array}{ccccccc}
    L(F)+lI_{n} & -1_{l}^{T}\otimes I_{n}& 0\\
   -1_{l}\otimes I_{n} & (L(H_{1})+(m-l+1)I_{l})\otimes I_{n}&-J_{l\times (m-l)}\otimes I_{n} \\
   0&-J_{(m-l)\times l}\otimes I_{n}&(L(H_{2})+lI_{m-l})
   \otimes I_{n}\\
  \end{array}
\right).
\end{array}
\]

We begin with the computation of $\{1\}$-inverse of $G=G[F;H_{v},l]$.

Let
$A=L(F)+lI_{n}$,
 $B=\left(
  \begin{array}{cccccccccccccccc}
   -1_{l}^{T}\otimes I_{n}& 0\\
 \end{array}
  \right)$, $B^{T}=\left(
  \begin{array}{cccccccccccccccc}
    -1_{l}\otimes I_{n} &\\
   0 &\\
 \end{array}
  \right)
$, and

\[
\begin{array}{crl}
D=\left(
  \begin{array}{cccccccccccccccc}
    (L(H_{1})+(m-l+1)I_{l})\otimes I_{n}&-J_{l\times (m-l)}\otimes I_{n} \\
   -J_{(m-l)\times l}\otimes I_{n}&(L(H_{2})+lI_{m-l})
   \otimes I_{n}\\
 \end{array}
  \right)
\end{array}.
\]

First we computer the $D^{-1}$. By Lemma 2.3, we have
\[
\begin{array}{crl}
 A_{1}-B_{1}D_{1}^{-1}C_{1}&=&
\left(L(H_{1})+(m-l+1)I_{l}\right)\otimes I_{n}-
 (-J_{l\times (m-l)}\otimes I_{n})((L(H_{2})+lI_{m-l})^{-1}
   \otimes I_{n})\\&&(-J_{(m-l)\times l}\otimes I_{n})\\
   &=&(L(H_{1})+(m-l+1)I_{l})\otimes I_{n}-1_{l}[1^{T}_{m-l}
   (L(H_{2})+lI_{m-l})^{-1}1_{m-l}]1^{T}_{l}\otimes I_{n}\\
   &=&(L(H_{1})+(m-l+1)I_{l})\otimes I_{n}-
   \frac{m-l}{l}J_{l\times l}\otimes I_{n}\\
  &=&[L(H_{1})+(m-l+1)I_{l}-\frac{m-l}{l}J_{l\times l}]
  \otimes I_{n},
\end{array}
\]
so $(A_{1}-B_{1}D_{1}^{-1}C_{1})^{-1}=[(L(H_{1})
+(m-l+1)I_{l}-\frac{m-l}{l}J_{l\times l}]^{-1}\otimes I_{n}$.

By Lemma 2.3, we have
\[
\begin{array}{crl}
 S^{-1}&=&(D_{1}-C_{1}A_{1}^{-1}B_{1})^{-1}\\
 &=&[\left(L(H_{2})+lI_{m-l}\right)\otimes I_{n}-
 (-J_{(m-l)\times l}\otimes I_{n})((L(H_{1})+(m-l+1)I_{l})^{-1}
   \otimes I_{n})\\&&(-J_{l\times (m-l)}\otimes I_{n})]^{-1}\\
   &=&[(L(H_{2})+lI_{m-l})\otimes I_{n}-(J_{(m-l)\times l}
   (L(H_{1})+(m-l+1)I_{l})^{-1}J_{l\times (m-l)})\otimes I_{n}]^{-1}\\
   &=&(L(H_{2})+lI_{m-l}-
   \frac{l}{m-l+1}J_{(m-l)\times (m-l)})^{-1}\otimes I_{n}.
\end{array}
\]

By Lemma 2.3, we have
\[
\begin{array}{crl}
 -A_{1}^{-1}B_{1}S^{-1}&=&-[(L(H_{1})+(m-l+1)I_{l})^{-1}\otimes I_{n}]
 (-J_{l\times(m-l)}\otimes I_{n})\\&&[L(H_{2})+lI_{m-l}-\frac{l}{m-l+1}J_{(m-l)\times (m-l)}]^{-1}\otimes I_{n}\\
 &=&\frac{1}{m-l+1}1_{l}\times\frac{m-l+1}{l}1^{T}_{m-l}\otimes I_{n}\\
 &=&\frac{1}{l}J_{l\times(m-l)}\otimes I_{n}.
\end{array}
\]

Similarly, $-S^{-1}C_{1}A_{1}^{-1}=(-A_{1}^{-1}B_{1}S^{-1})^{T}=
\frac{1}{l}J_{(m-l)\times l}\otimes I_{n}$.
So

$D^{-1}$
\begin{eqnarray*}
 &=&
 \left(
  \begin{array}{cccccccccccccccc}
   [(L(H_{1})+(m-l+1)I_{l}-\frac{m-l}{l}J_{l\times l}]^{-1}\otimes I_{n}
   &\frac{1}{l}J_{l\times(m-l)}\otimes I_{n}\\
   \frac{1}{l}J_{(m-l)\times l}\otimes I_{n}&[(L(H_{2})+lI_{m-l}-
   \frac{l}{m-l+1}J_{(m-l)\times (m-l)})]^{-1}\otimes I_{n}\\
 \end{array}
  \right).
\end{eqnarray*}

Now we compute the $\{1\}$-inverse of $G[F;H_{v},l]$. By Lemma 2.5, we have
\[
\begin{array}{crl}
 H&=&A-BD^{-1}B^{T}\\
 &=&L(F)+lI_{n}-\left(
  \begin{array}{cccccccccccccccc}
     -1_{l}^{T}\otimes I_{n} & 0\\
 \end{array}
  \right)\\&&\left(
  \begin{array}{cccccccccccccccc}
   [(L(H_{1})+(m-l+1)I_{l}-\frac{m-l}{l}J_{l\times l}]^{-1}\otimes I_{n}
   &\frac{1}{l}J_{l\times(m-l)}\otimes I_{n}\\
   \frac{1}{l}J_{(m-l)\times l}\otimes I_{n}&[L(H_{2})+lI_{m-l}-
   \frac{l}{m-l+1}J_{(m-l)\times (m-l)}]^{-1}\otimes I_{n}\\
 \end{array}
  \right)\\&&\left(
  \begin{array}{cccccccccccccccc}
     -1_{l}\otimes I_{n} \\
     0  \\
 \end{array}
  \right)\\
  &=&L(F)+lI_{n}-(1_{l}^{T}\otimes I_{n})[L(H_{1})+(m-l+1)I_{l}-
  \frac{m-l}{l}J_{l\times l}]^{-1}\otimes I_{n}](1_{l}\otimes I_{n})\\
 &=&L(F)+lI_{n}-lI_{n}\\
 &=&L(F),
\end{array}
\]
so $H^{\#}=L^{\#}(F)$.

According to Lemma 2.6, we calculate $-H^{\#}BD^{-1}$ and $-D^{-1}B^{T}H^{\#}$.

$-H^{\#}BD^{-1}$
\[
\begin{array}{crl}
&=&-L^{\#}(F)\left(
  \begin{array}{cccccccccccccccc}
     -1_{l}^{T}\otimes I_{n} & 0\\
 \end{array}
  \right)
\\&&\left(
  \begin{array}{cccccccccccccccc}
   [(L(H_{1})+(m-l+1)I_{l}-\frac{m-l}{l}J_{l\times l}]^{-1}\otimes I_{n}
   &\frac{1}{l}J_{l\times(m-l)}\otimes I_{n}\\
   \frac{1}{l}J_{(m-l)\times l}\otimes I_{n}&[(L(H_{2})+lI_{m-l}-
   \frac{l}{m-l+1}J_{(m-l)\times (m-l)}]^{-1}\otimes I_{n}\\
 \end{array}
  \right)\\
  &=&\left(L^{\#}(F)(1^{T}_{l}\otimes I_{n}),
  L^{\#}(F)(1^{T}_{m-l}\otimes I_{n})\right)
\end{array}
\]
and
\[
\begin{array}{crl}
-D^{-1}B^{T}H^{\#}
  &=&\left(
  \begin{array}{cccccccccccccccc}
 (1_{l}\otimes I_{n})L^{\#}(F)& \\
   (1_{m-l}\otimes I_{n})L^{\#}(F)&\\
 \end{array}
  \right).
\end{array}
\]

We are ready to compute the $D^{-1}B^{T}H^{\#}BD^{-1}$.
\[
\begin{array}{crl}
D^{-1}B^{T}H^{\#}BD^{-1}
&=&\left(\begin{array}{cccccccccccccccc}
 (1_{l}\otimes I_{n})L^{\#}(F)& \\
   (1_{m-l}\otimes I_{n})L^{\#}(F)& \\
 \end{array}\right)
\left((1^{T}_{l}\otimes I_{n}),
  (1^{T}_{m-l}\otimes I_{n})\right)\\
&=&\left(
  \begin{array}{cccccccccccccccc}
 (1_{l}\otimes I_{n})L^{\#}(F)(1_{l}^{T}\otimes I_{n})&
 (1_{l}\otimes I_{n})L^{\#}(F)(1_{m-l}^{T}\otimes I_{n})\\
(1_{m-l}\otimes I_{n})L^{\#}(F)(1^{T}_{l}\otimes I_{n})&(1_{m-l}\otimes I_{n})L^{\#}(F)(1^{T}_{m-l}\otimes I_{n}) \\
 \end{array}
  \right).\\
\end{array}\\
\]

Let $P=[(L(H_{1})+(m-l+1)I_{l}-\frac{m-l}{l}J_{l\times l}]\otimes I_{n}$, $Q=(L(H_{2})+lI_{m-l}-\frac{l}{m-l+1}J_{(m-l)\times (m-l)})^{-1}\otimes I_{n}$, $M=\frac{1}{l}
    J_{l\times (m-l)}\otimes I_{n}+(1_{l}\otimes I_{n})L^{\#}(F)(1^{T}_{m-l}\otimes I_{n})$, then based on
Lemma 2.6, the following matrix

$N=$
\begin{eqnarray}
\left(
  \begin{array}{cccccccccccccccc}
  L^{\#}(F)&L^{\#}(F)(1^{T}_{l}\otimes I_{n})&
  L^{\#}(F)(1^{T}_{m-l}\otimes I_{n})\\
    (1_{l}\otimes I_{n})L^{\#}(F)&P^{-1}+(1_{l}\otimes I_{n})L^{\#}(F)(1^{T}_{l}\otimes I_{n})&M\\
   (1_{m-l}\otimes I_{n})L^{\#}(F)&
   M^{T}&
    Q^{-1}+(1_{m-l}\otimes I_{n})L^{\#}(F)(1^{T}_{m-l}\otimes I_{n})\\
 \end{array}
  \right)
\end{eqnarray}
is a symmetric $\{1\}$-inverse of $G[F;H_{v},l]$.

For any $i,j\in V(F)$, by Lemma 2.1 and the Equation $(3.1)$, we have
\begin{eqnarray*}
r_{ij}(G[F;H_{v},l])=(L^{\#}(F))_{ii}+(L^{\#}(F))_{jj}
-2(L^{\#}(F))_{ij}=r_{ij}(F),
\end{eqnarray*}
as stated in (i).

For any $i\in V(F)$, $j\in V(H_{1})$, by Lemma 2.1 and the Equation $(3.1)$, we have
\begin{eqnarray*}
r_{ij}(G[F;H_{v},l])&=&(L^{\#}(F))_{ii}+
[(L(H_{1})+(m-l+1)I_{l}-\frac{m-l}{l}J_{l\times l})^{-1}\otimes I_{n}]_{jj}-2(L^{\#}(F))_{ij},
\end{eqnarray*}
as stated in (ii).

For any $i\in V(F)$, $j\in V(H_{2})$, by Lemma 2.1 and the Equation $(3.1)$, we have
\begin{eqnarray*}
r_{ij}(G[F;H_{v},l])&=&(L^{\#}(F))_{ii}+
[(L(H_{2})+lI_{m-l}-\frac{l}{m-l+1}J_{(m-l)\times (m-l)})^{-1}\otimes I_{n}]_{jj}-2(L^{\#}(F))_{ij},
\end{eqnarray*}
as stated in (iii).

For any $i\in V(H_{1})$, $j\in V(H_{2})$, by Lemma 2.1 and the Equation $(3.1)$, we have
\begin{eqnarray*}
r_{ij}(G[F;H_{v},l])&=&[(L(H_{1})+(m-l+1)I_{l}-\frac{m-l}{l}J_{l\times l})^{-1}\otimes I_{n}]_{ii}+\left[(L(H_{2})+lI_{m-l}\right.\\&&
\left.-\frac{l}{m-l+1}J_{(m-l)\times (m-l)})^{-1}\otimes I_{n}]_{jj}-2[(L(H_{1})+(m-l+1)I_{l}-\frac{m-l}{l}J_{l\times l})^{-1}\otimes I_{n}\right]_{ij},
\end{eqnarray*}
as stated in (iv).

For any $i\in V(H_{2})$, $j\in V(H_{1})$, by Lemma 2.1 and the Equation $(3.1)$, we have
\begin{eqnarray*}
r_{ij}(G[F;H_{v},l])&=&[(L(H_{2})+lI_{m-l}-\frac{l}{m-l+1}J_{(m-l)\times (m-l)})^{-1}\otimes I_{n}]_{ii}+\left[(L(H_{1})+(m-l+1)I_{l}\right.\\&&
\left.-\frac{m-l}{l}J_{l\times l})^{-1}\otimes I_{n}]_{jj}-2[(L(H_{2})+lI_{m-l}-\frac{l}{m-l+1}J_{(m-l)\times (m-l)})^{-1}\otimes I_{n}\right]_{ij},
\end{eqnarray*}
as stated in (v).

Now we compute the Kirchhoff index of $G[F;H_{v},l]$.
By Lemma 2.5, we have

$Kf(G[F;H_{v},l])$
\begin{eqnarray*}
&=&n(m+1)tr( N)-1^{T}N1\\
&=&n(m+1)\left[\left(tr(L^{\#}(F))+tr\left((L(H_{1})+(m-l+1)I_{l}+
\frac{m-l}{l}J_{l\times l})^{-1}\otimes I_{n})\right)+\right.\right.\\
&&\left.\left.+tr\left((L(H_{2})+lI_{m-l}+
\frac{l}{m-l+1}J_{(m-l)\times (m-l)})^{-1}\otimes I_{n}\right)\right)\right.+\\&&
\left.tr((1_{l}\otimes I_ {n})L^{\#}(F)(1^{T}_{l}\otimes I_{n}))
+tr((1_{m-l}\otimes I_{n})L^{\#}(F)(1^{T}_{m-l}\otimes I_{n}))\right]-1^{T}N1.
\end{eqnarray*}

Note that the eigenvalues of $(L(H_{2})+lI_{m-l})$
are $0+l, \nu_{2}(H_{2})+l,...,\nu_{m-l}(H_{2})+l$
and the eigenvalues of $J_{(m-l)\times(m-l)}$
is $(m-l),0^{(m-l-1)}$. Then
\begin{eqnarray*}
tr\left((L(H_{2})+lI_{m-l}-\frac{l}{m-l+1}J_{(m-l)\times
(m-l)})^{-1}\otimes I_{n}\right)
\end{eqnarray*}
\begin{eqnarray}
&=&n\sum_{i=2}^{m-l}\frac{1}{\nu_{i}(H_{2})+l}+\frac{nl}{m-l+1}.
\end{eqnarray}

Similarly,
$
tr\left((L(H_{1})+(m-l+1)I_{l}-
\frac{m-l}{l}J_{l\times l})^{-1}\otimes I_{n})\right)
=n\sum_{i=2}^{l}\frac{1}{\mu_{i}(H_{1})+(m-l+1)}+n.
$

It is easily obtained
\begin{eqnarray*}
tr\left((1_{l}\otimes I_ {n})L^{\#}(F)(1^{T}_{l}\otimes I_{n})\right)
+tr\left((1_{m-l}\otimes I_{n})L^{\#}(F)(1^{T}_{m-l}\otimes I_{n})\right)
\end{eqnarray*}
\begin{eqnarray*}
&=&tr\left(J_{l\times l}\otimes L^{\#}(F)\right)+
tr\left(J_{(m-l)\times (m-l)}\otimes L^{\#}(F)\right)
\end{eqnarray*}
\begin{eqnarray}
&=&ltr\left(L^{\#}(F)\right)+(m-l)tr\left(L^{\#}(F)\right)
=mtr\left(L^{\#}(F)\right).
\end{eqnarray}

Let $P=\left(L(H_{1})+(m-l+1)I_{l}-\frac{m-l}{l}J_{l\times l}\right)
\otimes I_{n}$, then
\begin{eqnarray*}
1^{T}P^{-1}1&=&\left(
  \begin{array}{cccccccccccccccc}
 1^{T}_{l}&1^{T}_{l}&\cdots &1^{T}_{l}\\
 \end{array}
  \right)
  \left(
  \begin{array}{cccccccccccccccc}
P^{-1}& 0 & 0&...&0\\
  0&P^{-1}&0&...&0\\
  0&0& ...&...&0\\
   0&0& 0&...&P^{-1}\\
 \end{array}
  \right)\left(
  \begin{array}{cccccccccccccccc}
   1_{l}\\
  1_{l}\\
  \cdots\\
  1_{l}\\
 \end{array}
  \right)\\
  \end{eqnarray*}
  \begin{eqnarray}
  &=&l1^{T}_{l}(L(H_{1})+(m-l+1)I_{l}-\frac{m-l}{l}J_{l\times l})^{-1}1_{l}=l^{2}
    \end{eqnarray}

Let $Q=\left(L(H_{2})+lI_{m-l}-\frac{l}{m-l+1}J_{(m-l)\times (m-l)}\right)$, then
\begin{eqnarray*}
1^{T}Q^{-1}1&=&\left(
  \begin{array}{cccccccccccccccc}
 1^{T}_{m-l}&1^{T}_{m-l}&\cdots &1^{T}_{m-l}\\
 \end{array}
  \right)
  \left(
  \begin{array}{cccccccccccccccc}
Q^{-1}& 0 & 0&...&0\\
  0&Q^{-1}&0&...&0\\
  0&0& ...&...&0\\
   0&0& 0&...&Q^{-1}\\
 \end{array}
  \right)\left(
  \begin{array}{cccccccccccccccc}
   1_{m-l}\\
  1_{m-l}\\
  \cdots\\
  1_{m-l}\\
 \end{array}
  \right)\\
  &=&(m-l)1^{T}_{m-l}(L(H_{2})+lI_{m-l}-\frac{l}{m-l+1}J_{(m-l)\times (m-l)})^{-1}1_{m-l}
  \end{eqnarray*}
 \begin{eqnarray}
  &=&(m-l)^{2}\frac{m-l+1}{l}.
  \end{eqnarray}
\begin{eqnarray*}
1^{T}_{ln}(1_{l}\otimes I_{n})L^{\#}(F)
(1^{T}_{l}\otimes I^{n})1_{ln}
  &=&\left(
  \begin{array}{cccccccccccccccc}
 1^{T}_{n}&1^{T}_{n}&\cdots&1^{T}_{n}\\
 \end{array}
 \right)
 \left(
 \begin{array}{cccccccccccccccc}
 I_{n}&\\
 I_{n}&\\
 \cdots&\\
 I_{n}&\\
 \end{array}
  \right)L^{\#}(F)\left(
  \begin{array}{cccccccccccccccc}
 I_{n}&I_{n}&\cdots&I_{n}\\
 \end{array}
 \right)
 \end{eqnarray*}
 \begin{eqnarray}
 \left(
 \begin{array}{cccccccccccccccc}
 1_{n}&\\
 1_{n}&\\
 \cdots&\\
 1_{n}&\\
 \end{array}
  \right)
  =n^{2}1^{T}_{n}L^{\#}(F)1_{n}=0,
\end{eqnarray}

Similarly,
$1^{T}((1_{l}\otimes I_{n})L^{\#}(F)(1^{T}_{m-l}\otimes I_{n})1=0,1^{T}((1_{m-l}\otimes I_{n})L^{\#}(F)
(1^{T}_{l}\otimes I_{n})1=0$
and $1^{T}((1_{m-l}\otimes I_{n})L^{\#}(F)
(1^{T}_{l}\otimes I_{n})1=0$

Plugging $(3.2)$,$(3.3)$,$(3.4)$, $(3.5)$ and (3.6) and the above equation
into $Kf(G[F;H_{v},l])$, we obtain the required result in (iv).

\section{Resistance distance and Kirchhoff index of $G[F,u_{1},u_{2},\ldots,u_{k};H_{v},l]$}
In this section, we consider the case when $F=F_{1}
\vee F_{2}$, where $F_{1}$ the subgraph of $F$ induced
by the vertices $u_{1},u_{2},\ldots,u_{k}$ and $F_{2}$
is the subgraph of $F$ induced by the vertices
$u_{k+1},u_{k+2},\ldots,u_{n}$. In this case, we give
the explicit formulate of resistance distance and Kirchhoff index of $G=G[F,u_{1},u_{2},\ldots,u_{k};H_{v},l]$.

{\bf Theorem 4.1} \ Let $G=G[F,u_{1},u_{2},\ldots,u_{k};H_{v},l]$
be the graph as described above. Let
$\sigma(F_{1})=(0=\alpha_{1},\alpha_{2},\ldots \alpha_{k})$, $\sigma(F_{2})=(0=\beta_{1},\beta_{2},\ldots \beta_{n-k})$ and
$\sigma(H_{1})=(0=\mu_{1},\mu_{2},\ldots \mu_{l})$, and
$\sigma(H_{2})=(0=\nu_{1},\nu_{2},\ldots \nu_{m-l})$.
Then $G$ have the resistance distance and Kirchhoff index
as follows:

(i)For any $i,j\in V(F_{1})$, we have
\begin{eqnarray*}
r_{ij}(G)&=&((L(F_{1})+(n-k)I_{k})^{-1}-
\frac{n-k}{k})_{ii}+((L(F_{1})+(n-k)I_{k})^{-1}-
\frac{n-k}{k})_{jj}\\&&
-2((L(F_{1})+(n-k)I_{k})^{-1}-
\frac{n-k}{k})_{ij}.
\end{eqnarray*}

(ii)For any $i,j\in V(F_{2})$, we have
\begin{eqnarray*}
r_{ij}(G)&=&(L(F_{2})+kI_{n-k})^{-1}_{ii}
+(L(F_{2})+kI_{n-k})^{-1}_{jj}
-2(L(F_{2})+kI_{n-k})^{-1}_{ij}.
\end{eqnarray*}

(iii)For any $i,j\in V(H_{1})$, we have
\begin{eqnarray*}
r_{ij}(G)&=&((L(H_{1})+(m-l+1)I_{l}-\frac{m-l}{l}J_{l\times l})\otimes I_{k})_{ii}+(L(H_{1})+(m-l+1)I_{l}-\frac{m-l}{l}J_{l\times l})\otimes I_{k})_{jj}\\&&
-2(L(H_{1})+(m-l+1)I_{l}-\frac{m-l}{l}J_{l\times l})\otimes I_{k})_{ij}.
\end{eqnarray*}

(iv)For any $i,j\in V(H_{2})$, we have
\begin{eqnarray*}
r_{ij}(G)&=&(L(H_{2})+lI_{m-l}-
   \frac{l}{m-l+1}J_{(m-l)\times (m-l)})\otimes I_{k})_{ii}+
(L(H_{2})+lI_{m-l}-\frac{l}{m-l+1}\\&&
  J_{(m-l)\times (m-l)})\otimes I_{k})_{jj}
-2(L(H_{2})+lI_{m-l}-
   \frac{l}{m-l+1}J_{(m-l)\times (m-l)})\otimes I_{k})_{ij}.
\end{eqnarray*}

(v)For any $i\in V(F)$,$j\in V(H_{1})$, we have
\begin{eqnarray*}
r_{ij}(G)&=&(L^{\#}(F))_{ii}+
[(L(H_{1})+(m-l+1)I_{l}-\frac{m-l}{l}j_{l\times l})]^{-1}\otimes I_{n}-2(L^{\#}(F))_{ij}.
\end{eqnarray*}

(vi)For any $i\in V(F)$, $j\in V(H_{2})$, we have
\begin{eqnarray*}
r_{ij}(G)&=&(L^{\#}(F))_{ii}+
[(L(H_{2})+lI_{m-l}-\frac{l}{m-l+1}j_{(m-l)\times (m-l)})]^{-1}\otimes I_{n}-2(L^{\#}(F))_{ij}.
\end{eqnarray*}

(vii)For any $i\in V(H_{1})$, $j\in V(H_{2})$, we have
\begin{eqnarray*}
r_{ij}(G)&=&[(L(H_{1})+(m-l+1)I_{l}-\frac{m-l}{l}j_{l\times l})^{-1}\otimes I_{n}]_{ii}+
[(L(H_{2})+lI_{m-l}-\\&&
\frac{l}{m-l+1}j_{(m-l)\times (m-l)})^{-1}\otimes I_{n}]_{jj}
-2[(L(H_{1})+(m-l+1)I_{l}-\frac{m-l}{l}j_{l\times l})^{-1}\otimes I_{n}]_{ij}.
\end{eqnarray*}

(viii)For any $i\in V(H_{2})$, $j\in V(H_{1})$, we have
\begin{eqnarray*}
r_{ij}(G)&=&[(L(H_{2})+lI_{m-l}-\frac{l}{m-l+1}j_{(m-l)\times (m-l)})^{-1}\otimes I_{n}]_{ii}+
[(L(H_{1})+(m-l+1)I_{l}-\\&&
\frac{m-l}{l}j_{l\times l})^{-1}\otimes I_{n}]_{jj}
-2[(L(H_{2})+lI_{m-l}-\frac{l}{m-l+1}j_{(m-l)\times (m-l)})^{-1}\otimes I_{n}]_{ij}.
\end{eqnarray*}

(ix)
\begin{eqnarray*}
Kf(G)&=&(n+mk)\left[
2\sum_{i=1}^{k}(\frac{1}{\alpha_{i}+(n-k)}-\frac{1}{(n-k)})
+\sum_{i=1}^{n-k}\frac{1}{\beta_{i}+k}+\right.\\&&\left.
(k\sum_{i=2}^{l}\frac{1}{\mu_{i}+(m-l+1)}+k)+
(k\sum_{i=2}^{m-l}\frac{1}{\nu_{i}+l}+\frac{l(2m-2l+1)}{m-l+1})
+k+\frac{k(m-l)}{l}\right]
-\\&&\left(l^{2}+\frac{(m-l)(m-l+1)}{l}+2k(m-l)\right).
\end{eqnarray*}

{\bf Proof} \ Let $v^{i}_{j}$ denote the $j$th vertex of $H_{v}$
in the $i$th copy of $H_{v}$ in $G$, for $i=1,2,\ldots,k$;
$j=1,2,\ldots,m$, and let $V_{j}(H_{v})=\{v^{1}_{j},v^{2}_{j},
\ldots,v^{k}_{j}\}$. Then $V(F_{1})\cup V(F_{2})\cup (\cup_{j=1}^{m}V_{j}(H_{v}))$
is a partition of the vertex set of $G=G[F,u_{1},u_{2},\ldots,u_{k};H_{v},l]$. Using this partition, the Laplacian
matrix of $G$ can be expressed as
\[
\begin{array}{crl}L(G)=
\left(
  \begin{array}{ccccccc}
    L_{1} & -J_{k\times(n-k)}& -1^{T}_{l}\otimes I_{k}& 0\\
   -J_{(n-k)\times k} & L_{2}&0&0 \\
   -1_{l}\otimes I_{k}&0&L_{3}&-J_{l\times(m-l)}\otimes I_{k}\\
   0&0&-J_{(m-l)\times l}\otimes I_{k}& L_{4}\\
  \end{array}
\right)
\end{array},
\]
where $L_{1}=L(F_{1})+(n-k+l)I_{k}$, $L_{2}=L(F_{2})+kI_{n-k}$,
$L_{3}=(L(H_{1})+(m-l+1)I_{l})\otimes I_{k}$, and
$L_{4}=(L(H_{2})+lI_{m-l})\otimes I_{k}.$

Let $A=L_{1}$,
 $B=\left(
  \begin{array}{ccccccc}
     -J_{k\times (n-k)}& -1^{T}_{l}\otimes I_{k}&0\\
  \end{array}
\right)$, $B^{T}=\left(
  \begin{array}{ccccccc}
   -J_{(n-k)\times k}&\\
   -1_{l}\otimes I_{k}&\\
   0&\\
  \end{array}
\right)$, and

\[
\begin{array}{crl}
 D=\left(
  \begin{array}{cccccccccccccccc}
     L_{2}&0&0 \\
   0&L_{3}&-J_{l\times(m-l)}\otimes I_{k}\\
   0&-J_{(m-l)\times l}\otimes I_{k}& L_{4}\\
 \end{array}
  \right)
\end{array}.
\]

First, we compute
$
D_{1}^{-1}=\left(
  \begin{array}{cccccccccccccccc}
   L_{3}&-J_{l\times(m-l)}\otimes I_{k}\\
   -J_{(m-l)\times l}\otimes I_{k}& L_{4}\\
 \end{array}
  \right)^{-1}
$.
By Lemma 2.3, we have
\[
\begin{array}{crl}
 A_{1}-B_{1}D_{1}^{-1}C_{1}
 &=&\left(L(H_{1})+(m-l+1)I_{l}\right)\otimes I_{k}-
 (-J_{l\times (m-l)}\otimes I_{k})((L(H_{2})+lI_{m-l})^{-1}
   \otimes I_{k})\\&&(-J_{(m-l)\times l}\otimes I_{k})\\
   &=&(L(H_{1})+(m-l+1)I_{l})\otimes I_{k}-1_{l}(1^{T}_{m-l}
   (L(H_{2})+lI_{m-l})^{-1}1_{m-l})1^{T}_{l}\otimes I_{k}\\
   &=&(L(H_{1})+(m-l+1)I_{l})\otimes I_{k}-
   \frac{m-l}{l}J_{l\times l}\otimes I_{k}\\
  &=&[L(H_{1})+(m-l+1)I_{l}-\frac{m-l}{l}J_{l\times l}]
  \otimes I_{k}
\end{array}
\]
so $(A_{1}-B_{1}D_{1}^{-1}C_{1})^{-1}=[(L(H_{1})+(m-l+1)I_{l}-\frac{m-l}{l}J_{l\times l}]^{-1}\otimes I_{k}$.

By Lemma 2.3, we have
\[
\begin{array}{crl}
 S^{-1}&=&(D_{1}-C_{1}A_{1}^{-1}B_{1})^{-1}\\
 &=&[\left(L(H_{2})+lI_{m-l}\right)\otimes I_{k}-
 (-J_{(m-l)\times l}\otimes I_{k})((L(H_{1})+(m-l+1)I_{l})^{-1}
   \otimes I_{k})\\&&(-J_{l\times (m-l)}\otimes I_{k}]^{-1}\\
   &=&[(L(H_{2})+lI_{m-l})\otimes I_{k}-(J_{(m-l)\times l}
   (L(H_{1})+(m-l+1)I_{l})^{-1}J_{l\times (m-l)})\otimes I_{k}]^{-1}\\
   &=&[L(H_{2})+lI_{m-l}-
   \frac{l}{m-l+1}J_{(m-l)\times (m-l)}]^{-1}\otimes I_{k}.
\end{array}
\]

By Lemma 2.3, we have
\[
\begin{array}{crl}
 -A_{1}^{-1}B_{1}S^{-1}&=&-[(L(H_{1})+(m-l+1)I_{l})^{-1}\otimes I_{k}]
 (-J_{l\times(m-l)}\otimes I_{k})\\&&[L(H_{2})+lI_{m-l}-\frac{l}{m-l+1}J_{(m-l)\times (m-l)}]^{-1}\otimes I_{k}\\
 &=&\frac{1}{m-l+1}1_{l}\times\frac{m-l+1}{l}1^{T}_{m-l}\otimes I_{k}\\
 &=&\frac{1}{l}J_{l\times (m-l)}\otimes I_{k}.
\end{array}
\]

Similarly, $-S^{-1}C_{1}A_{1}^{-1}=(-A_{1}^{-1}B_{1}S^{-1})^{T}=
\frac{1}{l}J_{(m-l)\times l}\otimes I_{k}$.
So

$D_{1}^{-1}=$
\begin{eqnarray*}
 \left(
  \begin{array}{cccccccccccccccc}
   [(L(H_{1})+(m-l+1)I_{l}-\frac{m-l}{l}J_{l\times l}]^{-1}\otimes I_{k}
   &\frac{1}{l}J_{l\times(m-l)}\otimes I_{k}\\
   \frac{1}{l}J_{(m-l)\times l}\otimes I_{k}&[(L(H_{2})+lI_{m-l}-
   \frac{l}{m-l+1}J_{(m-l)\times (m-l)})]^{-1}\otimes I_{k}\\
 \end{array}
  \right).
\end{eqnarray*}

Now we compute the $\{1\}$-inverse of $G[F,u_{1},u_{2},\ldots,u_{k};H_{v},l]$.
Let $P=[(L(H_{1})+(m-l+1)I_{l}-\frac{m-l}{l}J_{l\times l}]\otimes I_{k}$, $Q=[(L(H_{2})+lI_{m-l}-
   \frac{l}{m-l+1}J_{(m-l)\times (m-l)})]\otimes I_{k}$.
By Lemma 2.6, we have
\[
\begin{array}{crl}
H&=&A-BD^{-1}B^{T}\\
 &=&L(F_{1})+(n-k+l)I_{k}-\left(
  \begin{array}{cccccccccccccccc}
     -J_{k\times (n-k)}& -1^{T}_{l}\otimes I_{k}&0\\
 \end{array}
  \right)\\&&\left(
  \begin{array}{cccccccccccccccc}
  (L(F_{2})+kI_{n-k})^{-1}&0&0\\
   0&P^{-1}
   &\frac{1}{l}J_{l\times(m-l)}\otimes I_{k}\\
   0&\frac{1}{l}J_{(m-l)\times l}\otimes I_{k}&Q^{-1}\\
 \end{array}
  \right)\\&&\left(
  \begin{array}{cccccccccccccccc}
     -J_{(n-k)\times k} \\
     -1_{l}\otimes I_{k} \\
     0  \\
 \end{array}
  \right)\\
  &=&L(F_{1})+(n-k+l)I_{k}-\frac{n-k}{k}J_{k\times k}-lI_{k}\\
 &=&L(F_{1})+(n-k)I_{k}-\frac{n-k}{k}J_{k\times k},
\end{array}
\]
so $H^{\#}=(L(F_{1})+(n-k)I_{k}-\frac{n-k}{k}J_{k\times k})^{\#}$.
By Lemmma 2.4, $H^{\#}=(L(F_{1})+(n-k)I_{k})^{-1}-
\frac{1}{k(n-k)}J_{k\times k}$.

According to Lemma 2.6, we calculate $-H^{\#}BD^{-1}$ and $-D^{-1}B^{T}H^{\#}$.

$-H^{\#}BD^{-1}$
\[
\begin{array}{crl}
&=&-H^{\#}\left(
  \begin{array}{cccccccccccccccc}
    -J_{k\times (n-k)}& -1^{T}_{l}\otimes I_{k}&0\\
 \end{array}
  \right)
\\&&\left(
  \begin{array}{cccccccccccccccc}
   (L(F_{2})+kI_{n-k})^{-1}&0&0\\
   0&P^{-1}
   &\frac{1}{l}J_{l\times(m-l)}\otimes I_{k}\\
   0&\frac{1}{l}J_{(m-l)\times l}\otimes I_{k}&Q^{-1}\\
 \end{array}
  \right)\\
  &=&\left(
  \begin{array}{cccccccccccccccc}
   \frac{1}{k}H^{\#}J_{k\times(n-k)}
   &H^{\#}(1^{T}_{l}\otimes I_{k})&
   H^{\#}(1^{T}_{m-l}\otimes I_{k})\\
 \end{array}
  \right)
\end{array}
\]
and
\[
\begin{array}{crl}
-D^{-1}B^{T}H^{\#}
  &=&\left(
  \begin{array}{cccccccccccccccc}
 \frac{1}{k}J_{(n-k)\times k}H^{\#} & \\
   (1_{l}\otimes I_{k})H^{\#}&\\
   (1_{m-l}\otimes I_{k})H^{\#}\\
 \end{array}
  \right).\\
\end{array}
\]

We are ready to compute the $D^{-1}B^{T}H^{\#}BD^{-1}$.

$D^{-1}B^{T}H^{\#}BD^{-1}$
\[
\begin{array}{crl}
&=&\left(\begin{array}{cccccccccccccccc}
 \frac{1}{k}J_{(n-k)\times k}H^{\#} & \\
   (1_{l}\otimes I_{k})H^{\#}&\\
  (1_{m-l}\otimes I_{k})H^{\#}\\
 \end{array}\right)
\left(
  \begin{array}{cccccccccccccccc}
  \frac{1}{k}J_{k\times(n-k)}
   &(1^{T}_{l}\otimes I_{k})&
   (1^{T}_{m-l}\otimes I_{k})\\
 \end{array}
  \right)\\
&=&\left(
  \begin{array}{cccccccccccccccc}
 \frac{1}{k^{2}}J_{(n-k)\times k}H^{\#}J_{k\times (n-k)}&
 \frac{1}{k}J_{(n-k)\times k}H^{\#}(1^{T}_{l}\otimes I_{k})&
 \frac{1}{k}J_{(n-k)\times k}H^{\#}(1^{T}_{m-l}\otimes I_{k})\\
\frac{1}{k}(1_{l}\otimes I_{k})H^{\#}J_{k\times (n-k)}&(1_{l}\otimes I_{k})H^{\#}(1^{T}_{l}\otimes I_{k})
&(1_{l}\otimes I_{k})H^{\#}(1^{T}_{m-l}\otimes I_{k}) \\
\frac{1}{k}(1_{m-l}\otimes I_{k})H^{\#}J_{k\times (n-k)}&(1_{m-l}\otimes I_{k})H^{\#}(1^{T}_{l}\otimes I_{k})&
(1_{m-l}\otimes I_{k})H^{\#}(1^{T}_{m-l}\otimes I_{k})\\
 \end{array}
  \right).\\
\end{array}\\
\]

Let $M=1^{T}_{m-l}\otimes I_{k}$, $N=1^{T}_{l}\otimes I_{k}$.
Based on Lemma 2.6, the following matrix
\begin{eqnarray}
\left(
  \begin{array}{cccccccccccccccc}
 H^{\#}& \frac{1}{k}H^{\#}J_{k\times(n-k)}&H^{\#}N&H^{\#}M\\
   \frac{1}{k}J_{(n-k)\times k}H^{\#}  & (L(F_{2})+kI_{n-k})^{-1}& 0&0\\
    N^{T}H^{\#}&0&P^{-1}+N^{T}H^{\#}N&N^{T}H^{\#}M+
    \frac{1}{l}J_{l\times(m-l)}\otimes I_{k}\\
    M^{T}H^{\#}&0&M^{T}H^{\#}N+\frac{1}{l}J_{(m-l)\times l}\otimes I_{k}&Q^{-1}+M^{T}H^{\#}M)\\
 \end{array}
  \right)
\end{eqnarray}
is a symmetric $\{1\}$-inverse of $G=G[F,u_{1},u_{2},\ldots,u_{k};H_{v},l]$,
$P=[(L(H_{1})+(m-l+1)I_{l}-\frac{m-l}{l}J_{l\times l}]\otimes I_{k}$, $Q=[(L(H_{2})+lI_{m-l}-
   \frac{l}{m-l+1}J_{(m-l)\times (m-l)})]\otimes I_{k}$.

For any $i,j\in V(F_{1})$, by Lemma 2.1 and the Equation $(4.7)$, we have
\begin{eqnarray*}
r_{ij}(G)&=&((L(F_{1})+(n-k)I_{k})^{-1}-
\frac{1}{k(n-k)})_{ii}+((L(F_{1})+(n-k)I_{k})^{-1}-
\frac{1}{k(n-k)})_{jj}\\&&
-2((L(F_{1})+(n-k)I_{k})^{-1}-
\frac{1}{k(n-k)})_{ij},
\end{eqnarray*}
as stated in (i).

For any $i,j\in V(F_{2})$, by Lemma 2.1 and the Equation $(4.7)$, we have
\begin{eqnarray*}
r_{ij}(G)&=&(L(F_{2})+kI_{n-k})^{-1}_{ii}
+(L(F_{2})+kI_{n-k})^{-1}_{jj}
-2(L(F_{2})+kI_{n-k})^{-1}_{ij},
\end{eqnarray*}
as stated in (ii).

For any $i,j\in V(H_{1})$, by Lemma 2.1 and the Equation $(4.7)$, we have
\begin{eqnarray*}
r_{ij}(G)&=&((L(H_{1})+(m-l+1)I_{l}-\frac{m-l}{l}J_{l\times l})\otimes I_{k})^{-1}_{ii}+(L(H_{1})+(m-l+1)I_{l}-\frac{m-l}{l}J_{l\times l})\otimes I_{k})^{-1}_{jj}\\&&
-2(L(H_{1})+(m-l+1)I_{l}-\frac{m-l}{l}J_{l\times l})\otimes I_{k})^{-1}_{ij},
\end{eqnarray*}
as stated in (iii).

For any $i,j\in V(H_{2})$, by Lemma 2.1 and the Equation $(4.7)$, we have
\begin{eqnarray*}
r_{ij}(G)&=&(L(H_{2})+lI_{m-l}-
   \frac{l}{m-l+1}J_{(m-l)\times (m-l)})\otimes I_{k})^{-1}_{ii}+
(L(H_{2})+lI_{m-l}-\frac{l}{m-l+1}\\&&
  J_{(m-l)\times (m-l)})\otimes I_{k})^{-1}_{jj}
-2(L(H_{2})+lI_{m-l}-
   \frac{l}{m-l+1}J_{(m-l)\times (m-l)})\otimes I_{k})^{-1}_{ij},
\end{eqnarray*}
as stated in (iv).

For any $i\in V(F)$,$j\in V(H_{1})$ by Lemma 2.1 and the Equation $(4.7)$, we have
\begin{eqnarray*}
r_{ij}(G)&=&(L^{\#}(F))_{ii}+
[(L(H_{1})+(m-l+1)I_{l}-\frac{m-l}{l}j_{l\times l})^{-1}\otimes I_{n}]_{jj}-2(L^{\#}(F))_{ij},
\end{eqnarray*}
as stated in (v).

For any $i\in V(F)$, $j\in V(H_{2})$, by Lemma 2.1 and the Equation $(4.7)$, we have
\begin{eqnarray*}
r_{ij}(G)&=&(L^{\#}(F))_{ii}+
[(L(H_{2})+lI_{m-l}-\frac{l}{m-l+1}j_{(m-l)\times (m-l)})^{-1}\otimes I_{n}]_{jj}-2(L^{\#}(F))_{ij},
\end{eqnarray*}
as stated in (vi).

For any $i\in V(H_{1})$, $j\in V(H_{2})$, by Lemma 2.1 and the Equation $(4.7)$, we have
\begin{eqnarray*}
r_{ij}(G)&=&[(L(H_{1})+(m-l+1)I_{l}-\frac{m-l}{l}j_{l\times l})^{-1}\otimes I_{n}]_{ii}+
[(L(H_{2})+lI_{m-l}-\\&&
\frac{l}{m-l+1}j_{(m-l)\times (m-l)})^{-1}\otimes I_{n}]_{jj}
-2[(L(H_{1})+(m-l+1)I_{l}-\frac{m-l}{l}j_{l\times l})^{-1}\otimes I_{n}]_{ij},
\end{eqnarray*}
as stated in (vii).

For any $i\in V(H_{2})$, $j\in V(H_{1})$, by Lemma 2.1 and the Equation $(4.7)$, we have
\begin{eqnarray*}
r_{ij}(G)&=&[(L(H_{2})+lI_{m-l}-\frac{l}{m-l+1}j_{(m-l)\times (m-l)})^{-1}\otimes I_{n}]_{ii}+
[(L(H_{1})+(m-l+1)I_{l}-\\&&
\frac{m-l}{l}j_{l\times l})^{-1}\otimes I_{n}]_{jj}
-2[(L(H_{2})+lI_{m-l}-\frac{l}{m-l+1}j_{(m-l)\times (m-l)})^{-1}\otimes I_{n}]_{ij},
\end{eqnarray*}
as stated in (viii).

Now we compute the Kirchhoff index of $G[F,u_{1},u_{2},\ldots,u_{k};H_{v},l]$.

$Kf(G[F,u_{1},u_{2},\ldots,u_{k};H_{v},l])$
\begin{eqnarray*}
&=&(n+mk)tr( N)-1^{T}N1\\
&=&(n+mk)\left(tr((L(F_{1})+(n-k)I_{k})^{-1}-
\frac{1}{k(n-k)}J_{k\times k})+tr(L(F_{2})+kI_{n-k})^{-1}+\right.\\
&&\left.+ktr(L(H_{1})+(m-l+1)I_{l}-\frac{m-l}{l}J_{l\times l})^{-1}+ktr(L(H_{2})+lI_{m-l}-
   \frac{l}{m-l+1}J_{(m-l)\times (m-l)})^{-1}\right.\\&&
\left.+\frac{1}{l}tr(J_{l\times(m-l)}\otimes I_{k})
+\frac{1}{l}tr(J_{(m-l)\times l}\otimes I_{k})+tr(N^{T}H^{\#}N)
+tr(M^{T}H^{\#}M)\right)-1^{T}N1.
\end{eqnarray*}

Note that the eigenvalues of $(L(F_{1})+(n-k)I_{k})$
are $\alpha_{1}+(n-k), \alpha_{2}+(n-k),...,\alpha_{k}+(n-k)$. Then
\begin{eqnarray}
tr((L(F_{1})+(n-k)I_{k})^{-1}-
\frac{1}{k(n-k)}J_{k\times k})=\sum_{i=1}^{k}\frac{1}{\alpha_{i}+(n-k)}-\frac{k}{k(n-k)}.
\end{eqnarray}

Similarly,
$
tr\left((L(F_{2})+kI_{n-k})^{-1}\right)
=\sum_{i=1}^{n-k}\frac{1}{\beta_{i}+k}.
$

Note that the eigenvalues of $(L(H_{1})+(m-l+1)I_{l}$
are $0+(m-l+1), \mu_{2}(H_{1})+(m-l+1),...,\mu_{l}(H_{1})+(m-l+1)$
and the eigenvalues of $J_{(m-l)\times(m-l)}$
is $(m-l),0^{(m-l-1)}$. Then
\begin{eqnarray*}
tr\left(L(H_{1})+(m-l+1)I_{l}+
\frac{m-l}{l}J_{l\times l})^{-1}\otimes I_{k}\right)^{-1}
\end{eqnarray*}
\begin{eqnarray}
&=&k\sum_{i=2}^{l}\frac{1}{\mu_{i}+(m-l+1)}+k.
\end{eqnarray}

Similarly,
\begin{eqnarray}
tr\left((L(H_{2})+lI_{m-l}-\frac{l}{m-l+1}J_{(m-l)\times
(m-l)})\otimes I_{k}\right)^{-1}
=k\sum_{i=2}^{m-l}\frac{1}{\nu_{i}+l}+
\frac{kl(2m-2l+1)}{m-l+1}
\end{eqnarray}

It is easily obtained that
$tr(J_{l\times(m-l)}\otimes I_{k})=lk,
tr(J_{(m-l)\times l}\otimes I_{k})=(m-l)k,$
and
$
tr(N^{T}H^{\#}N)+tr(M^{T}H^{\#}M)
=tr(J_{l\times l}\otimes H^{\#})+tr(J_{(m-l)\times (m-l)}\otimes H^{\#})=ltr(H^{\#})+(m-l)tr(H^{\#})=mtr(H^{\#})$.

Since $1^{T}_{k}H^{\#}=1^{T}_{k}[(L(F_{1})+(n-k)I_{k})^{-1}-
\frac{1}{k(n-k)}J_{k\times k}]=\frac{1}{n-k}1^{T}_{k}-
\frac{1}{k(n-k)}1^{T}_{k}=0$, then
\begin{eqnarray*}
1^{T}N1&=&1^{T}(L(F_{2})+kI_{n-k})^{-1}1+1^{T}P^{-1}1+1^{T}Q^{-1}1
+1^{T}N^{T}H^{\#}N1+1^{T}N^{T}H^{\#}M1\\&&
+1^{T}M^{T}H^{\#}N1+1^{T}M^{T}H^{\#}M1
+\frac{1}{l}1^{T}(J_{l\times(m-l)}\otimes I_{k})1
+\frac{1}{l}1^{T}(J_{(m-l)\times l}\otimes I_{k})1.
\end{eqnarray*}

By the process of Theorem 3.1, we have
\begin{eqnarray*}
1^{T}P^{-1}1=l^{2},1^{T}Q^{-1}1=
(m-l)\frac{m-l+1}{l}.
\end{eqnarray*}
\[
\begin{array}{crl}
1^{T}(N^{T}H^{\#}N)1
&=&1^{T}_{lk}(1_{l}\otimes I_{k})H^{\#}
(1^{T}_{l}\otimes I^{k})1_{lk}\\
  &=&\left(
  \begin{array}{cccccccccccccccc}
 1^{T}_{k}&1^{T}_{k}&\cdots&1^{T}_{k}\\
 \end{array}
 \right)
 \left(
 \begin{array}{cccccccccccccccc}
 I_{k}&\\
 I_{k}&\\
 \cdots&\\
 I_{k}&\\
 \end{array}
  \right)H^{\#}\left(
  \begin{array}{cccccccccccccccc}
 I_{k}&I_{k}&\cdots&I_{k}\\
 \end{array}
 \right)\\&&
 \left(
 \begin{array}{cccccccccccccccc}
 1_{k}&\\
 1_{k}&\\
 \cdots&\\
 1_{k}&\\
 \end{array}
  \right)
  =k^{2}1^{T}_{k}H^{\#}1_{k}=0,
\end{array}
\]

Similarly, $1^{T}(M^{T}H^{\#}M)1=0$, $1^{T}N^{T}H^{\#}M1=0$
and $1^{T}M^{T}H^{\#}N1=0$.
\[
\begin{array}{crl}
1^{T}(J_{l\times(m-l)}\otimes I_{k})1
  &=&\left(
  \begin{array}{cccccccccccccccc}
 1^{T}_{k}&1^{T}_{k}&\cdots&1^{T}_{k}\\
 \end{array}
 \right)
 \left(
 \begin{array}{cccccccccccccccc}
 I_{k}&I_{k}&\cdots I_{k}&\\
 I_{k}&I_{k}&\cdots I_{k}&\\
 \cdots&\cdots&\\
 I_{k}&I_{k}&\cdots I_{k}&\\
 \end{array}
  \right)
 \left(
 \begin{array}{cccccccccccccccc}
 1_{k}&\\
 1_{k}&\\
 \cdots&\\
 1_{k}&\\
 \end{array}
  \right)\\
  &=&lk(m-l).
\end{array}
\]

Similarly, $1^{T}(J_{(m-l)\times l}\otimes I_{k})
=lk(m-l)$.

Plugging the above equation into $Kf(G[F,u_{1},u_{2},\ldots,u_{k};H_{v},l])$,
we obtain the required result in (iv).

\section{Conclusion}
In this paper,
using the Laplacian generalized inverse
approach, we obtained the resistance distance
and Kirchhoff indices of $G= G[F,u_{1},...,u_{k}, H_{v}]$ in terms of the resistance distance and Kirchhoff index
$F$ and $H_{uv}$, respectively. 

This article has been reviewed in Applied Mathematics-A Journal of Chinese Universities on September 26, 2018.


\noindent{\bf Acknowledgment:} 
This work was
supported by the National Natural Science Foundation of China
(Nos. 11461020),  the Research Foundation of the Higher
Education Institutions of Gansu Province, China (2018A-093), the Science and Technology Plan of Gansu Province (18JR3RG206) and
Research and Innovation Fund Project of President of Hexi University (XZZD2018003).


\begin{thebibliography}{99}

\bibitem{KR}D. J. Klein, M. Randi$\acute{c}$, Resistance distance, \emph{J. Math. Chem}, 12(1993)81-95.

\bibitem{JZ}S.B. Huang, J. Zhou and C.J. Bu, Some results on Kirchhoff index and
degree-Kirchhoff index, \emph{MATCH Commun. Math. Comput. Chem.} 75(2016)207-222.



\bibitem{JB1}J. D. Cao, J. B. Liu and S. H. Wang, Resistance distances in corona and neighborhood corona networks based on Laplacian generalized inverse approach
, \emph{J. Algebra. Appl.}(2019),1950053.

\bibitem{JB2}J. B. Liu, X. F. Pan, L. Yu and D. Li, Complete
characterization of bicyclic
graphs with minimal Kirchhoff index, \emph{Discrete. Appl. Math.} 200(2016)95-107.

\bibitem{SWZhB}L. Z. Sun, W. Z. Wang,  J. Zhou and C. J. Bu, Some results on resistance distances and resistance matrices, \emph{Linear and Multilinear Algebra}, 63(3)(2015)523-533.

\bibitem{R. B. Bapat}R. B. Bapat, \emph{Graphs and matrices}, Universitext, Springer/Hindustan Book Agency, London/New Delhi, 2010.

 \bibitem{ChZh}H. Y. Cheng and F. J. Zhang, Resistance distance and the normalized Laplacian spectrum. \emph{Discrete Appl. Math}, 155(2007)654-661.

 \bibitem{XG}W. J. Xiao and I. Gutman, Resistance distance and Laplacian spectrum, \emph{Theor. Chem. Acc},
110(2003), 284-289.

\bibitem{YK}Y. J. Yang and D. J. Klein, A recursion formula for resistance distances and its applications, \emph{Discrete Appl. Math}, 161(2013)2702-2715.

\bibitem{YK}Y. J. Yang and D. J. Klein, Resistance distance-based graph invariants of subdivisions and triangulations of graphs, \emph{Discrete Appl. Math}, 181(2015)260-274.


\bibitem{SB}S. Barik, On the Laplacian spectra of graphs with pockets,
\emph{Linear and Multilinear Algebra}, 56(2008), 481-490.

\bibitem{SB2018}S. Barik and G. Sahoo, Some results on the Laplacian spectra of graphs with pockets,
\emph{Electronic Notes in Discrete Mathematics}, 63(2017)219-228.

\bibitem{BG}R. B. Bapat and S. Gupta, Resistance distance in wheels and fans, \emph{Indian J. Pure Appl.Math}, 41(2010)1-13.


\bibitem{BYZhZH}C. J. Bu, Resistance distance in subdivision-vertex join and subdivision-edge join of graphs, \emph{Linear Algebra Appl.}, 458(2014)454-462.


 \bibitem{LZhB}X. G. Liu, J. Zhou and C. J. Bu, Resistance distance and Kirchhoff index of $R$-vertex join and $R$-edge join of two graphs, \emph{Discrete Appl. Math}, 187(2015) 130-139.



 \bibitem{A. Ben-Israel}A. Ben-Israel, T. N. E. Greville, \emph{Generalized inverses: theory and applications}, 2nd ed., Springer, New York, 2003.

\bibitem{BSZHW}C. J. Bu, L. Z. Sun, J. Zhou and Y. M. Wei, A note on block representations of the group inverse of Laplacian matrices, \emph{Electron. J. Linear Algebra}, 23(2012) 866-876.


\bibitem{FZZh}F. Z. Zhang, \emph{The Schur Complement and Its Applications}, Springer-Verlag, New York, 2005.



\bibitem{QL}Qun Liu, Some results of resistance distance and Kirchhoff index of vertex-edge corona for graphs, \emph{Advances in Mathematics(China)}, 45(2)(2016)176-183.










\end{thebibliography}
\end{document}